# Kalmár-style constructive completeness proofs for classical positive propositional calculi

Luiz Henrique Lopes dos Santos

USP/CNPq

## Preliminary remark

In the 1970s, concurrently with my doctoral research on Frege's philosophy of logic, I did some research in mathematical logic, under the supervision of Newton da Costa. At that time, Andrea Loparic and I had very frequent personal and professional contact, facilitated by our very close friendship.

As a routine, we discussed and critically reviewed each other's logical work in progress, as well as the drafts of papers intended to expose its results. In 1979, Andrea participated intensely in the elaboration process of the last mathematical logic paper I wrote, which dealt with the structure of some positive propositional calculi and their mutual connections.

This paper was published the following year in the Proceedings of the Third Brazilian Conference on Mathematical Logic. Unfortunately, several typographical errors, particularly those concerning logical symbols, made this first printed version of the paper practically unreadable. As a tribute to Andrea, I think it is worth publishing now a corrected (and slightly modified) version of it - which as a matter of fact will be its first truly intelligible printed version.

## 1. Introduction

The completeness of the axiomatic systems usually presented as formal descriptions of classical positive propositional logics (propositional logics without negation) is known to be ensured through distinct proof procedures. In Rasiowa 1974, for instance, we find several examples of how this can be done through algebraic methods. Pollock 1971 and Schumm 1975 display variants of Henkin's completeness proof procedures especially appropriate for applications to positive calculi.

In this paper, we shall settle the completeness of some classical positive propositional calculi (positive propositional calculi in which the so-called Peirce's law holds) by resorting to a close adaptation of Kalmár's completeness proof procedure, which was presented and applied to the classical propositional calculus with negation in Kalmár 1934-1935, and also in Mendelson 1964, for example.

First of all, we shall employ this adaptation to establish the completeness of the most familiar axiomatic characterization of the classical logic of material implication and disjunction. Next, we shall demonstrate that the completeness of the so-called classical positive propositional calculus (the classical calculus of material implication, disjunction, and conjunction) can be proved either through essentially the same adaptation of Kalmár's procedure or as an immediate consequence of the conjunction of our first result and the circumstance of each positive formula being reducible to a syntactically and semantically equivalent one exhibiting a particular normal form.

Finally, by defining disjunction in terms of material implication and so mirroring the classical calculus of material implication and disjunction in the classical implicative calculus, the completeness of the latter, as well as the completeness of the classical calculus of material implication and conjunction, will be shown to follow from the completeness of the former in quite simple and direct ways.

It is worth noting that all these completeness proofs will be constructive, that is to say, we can easily extract from each of them an effective method for producing formal deductions of all theses of the axiomatic system concerned.

## 2. Languages, calculi, and definitions

Throughout this paper, we shall treat *languages* as sets of expressions, to be called their *formulas*. Let L be the language whose formulas are denumerably many *atomic formulas* $p_1,\ldots,p_n,\ldots$, as well as what results from substituting formulas for the letters A and B in the schemes $(A \rightarrow B)$, $(A \vee B)$, and $(A \& B)$. The symbols $\rightarrow$, $\vee$ and $\&$ are the *connectives* of material implication (to be named from now on simply as implication), disjunction, and conjunction respectively. Capital letters A, B, C, D, and E, with or without numerical subscripts, will be used as syntactical variables for formulas of L unless restrictions are explicitly stated. The (proper and non-proper) parts of a formula will be called its subformulas. The formulas of L are supposed to be arranged in a fixed decidable linear order R (using Gödel numbering, for example).

For the sake of brevity, we may omit the outermost parentheses of an isolated formula. Sequences of formulas joined by the same connective are to be read as though they had the parentheses that could be written according to the familiar rule of association to the right. The schematic indication of such sequences by expressions of the kind of $A_1 \vee \ldots \vee A_n$, for example, is to be understood as comprising the case $n = 1$; in this case, $A_1 \vee \ldots \vee A_n$ will indicate the same as $A_1$. Similarly, the schematic expression for a sequence of formulas $A_1,\ldots,A_n$ will indicate the same as $A_1$, if $n = 1$.

The connective of equivalence is not a symbol of L. It will be used only for the sake of brevity, according to the usual definition:

$$A \leftrightarrow B =_{df} (A \rightarrow B) \,\&\, (B \rightarrow A).$$

A formula of L that contains no occurrence of a connective other than $\rightarrow$ will be called an *implicative formula*; one that contains no occurrence of $\vee$ will be called an *implicative-conjunctive formula*; one that contains no occurrence of & will be called an *implicative-disjunctive formula*.

The language of the *classical positive propositional calculus* is L and its axiom schemes and inference rule are:

(Ax1) $A \rightarrow B \rightarrow A$

(Ax2) $(A \rightarrow B \rightarrow C) \rightarrow (A \rightarrow B) \rightarrow A \rightarrow C$

(Ax3) $((A \rightarrow B) \rightarrow A) \rightarrow A$ (Peirce's law)

(Ax4) $A \rightarrow (A \vee B)$

(Ax5) $A \rightarrow (B \vee A)$

(Ax6) $(A \rightarrow C) \rightarrow (B \rightarrow C) \rightarrow (A \vee B) \rightarrow C$

(Ax7) $(A \,\&\, B) \rightarrow A$

(Ax8) $(A \,\&\, B) \rightarrow B$

(Ax9) $A \to B \to (A \& B)$

(MP) infer B from $A \to B$ and A.

The language of the *classical implicative propositional calculus* is the set of implicative formulas of L and its axiom schemes and inference rule are (Ax1)-(Ax3) and MP. The language of the *classical implicative-disjunctive propositional calculus* is the set of implicative-disjunctive formulas of L and its axiom schemes and inference rule are (Ax1)-(Ax6) and MP. The language of the *classical implicative-conjunctive propositional calculus* is the set of implicative-conjunctive formulas of L and its axiom schemes and inference rule are (Ax1)-(Ax3), (Ax7)-(Ax9), and MP.

The classical positive propositional calculus is clearly an extension of the classical implicative-disjunctive propositional calculus and also of the classical implicative-conjunctive propositional calculus. Each of these calculi is clearly an extension of the classical implicative propositional calculus.

The notation $\vdash_i A$ will mean that A is a thesis of the classical implicative propositional calculus. If K is a (possibly empty) set of implicative formulas, then the notation $K \vdash_i A$ will mean that A is formally derivable from the members of K in this calculus; the notation $B_1,\ldots,B_n \vdash_i A$ will be used instead of $\{B_1,\ldots,B_n\} \vdash_i A$. The signs $\vdash_{id}$, $\vdash_{ic}$ and $\vdash$ will be similarly used concerning the classical implicative-disjunctive propositional calculus, the classical implicative-conjunctive propositional calculus, and the classical positive propositional calculus respectively. Formulas A and B of L will be called *syntactical equivalents* in a calculus if and only if A is formally derivable from {B} and B is formally derivable from {A} in this calculus.

A *truth value assignment* for L is a function that maps L into the set {T, F} of truth values and fulfills the conditions classically associated with the connectives of implication, disjunction, and conjunction. The capital letter V will be used as a variable ranging over truth value assignments for L unless restrictions are explicitly stated. We define the concepts of tautology and tautological consequence as usual. The notation $\models A$ means that A is a tautology and $B_1,\ldots,B_n \models A$ means that A is a tautological consequence of $\{B_1,\ldots,B_n\}$. Formulas A and B of L will be called *semantical equivalents* if and only if $\models A \leftrightarrow B$.

Let A be a formula of L and let V be a truth value assignment for L. The notation $\Gamma$ [V; A] will refer to the set of atomic subformulas B of A such that V(B) = T; the notation $\Delta$ [V; A] will refer to the set of atomic subformulas B of A such that V(B) = F.

Let K be a finite set of formulas of L; if K is non-empty, let $B_1,\ldots,B_n$ be all the distinct elements of K arranged according to their relative positions in the linear order R. We shall denote by $(K)^A$ the formula A itself, if K is empty, and the formula $B_1 \vee \ldots \vee B_n \vee A$ otherwise; and we shall denote by $(K)^{\sim A}$ the formula A itself, if K is empty, and the formula $A \to (B_1 \vee \ldots \vee B_n)$ otherwise.

The schematic propositions 2.1-2.27 below can be easily proven. It is clear that schematic propositions 2.5-2.10 also hold for the classical implicative-disjunctive, classical implicative-conjunctive, and classical positive propositional calculi, and schematic propositions 2.11-2.24 also hold for the classical positive propositional calculus.

2.1. If B is a subformula of A, then $\Gamma$ (V, B) $\subseteq$ $\Gamma$ (V; A) and $\Delta$ (V, B) $\subseteq$ $\Delta$ (V; A).

2.2. If V(A) = F, then there is an atomic subformula B of A such that V(B) = F.

2.3. Deduction Theorem (DT) for all the propositional calculi specified above.

2.4. Semantical and syntactical rules of substitution of equivalents for all the propositional calculi specified above (RSE).

2.5. $\vdash_i A \to A$.

2.6. $A \to B, B \to C \vdash_i A \to C$.

2.7. $\vdash_i A \to (A \to B) \to B$.

2.8. $\vdash_i A \to (B \to A) \to A$.

2.9. $\vdash_i (A \to (A \to B)) \to A \to B$.

2.10. $(A \to B) \to B, (A \to C) \vdash_i (C \to B) \to B$.

2.11. $\vdash_{id} A \vee (A \to B)$.

2.12. $A \vee B, A \to C, B \to D \vdash_{id} C \vee D$.

2.13. $A \to B \vdash_{id} B \vee (A \to C)$

2.14. $A \vdash_{id} B \vee (C \to A)$.

2.15. $(A_1 \vee \ldots \vee A_n) \vee B$ and $A_1 \vee \ldots \vee A_n \vee B$ are syntactical equivalents in the classical implicative-disjunctive propositional calculus.

2.16. $A_1 \vee \ldots \vee A_n \vdash_{id} B_1 \vee \ldots \vee B_k$, if each of the $A_1, \ldots, A_n$ is one of the $B_1, \ldots, B_k$.

2.17. $A \vee B, C \to A \vdash_{id} (B \to C) \to A$.

2.18. $A \vee B, A \to B \vdash_{id} B$.

2.19. $(A \to B) \to B \vdash_{id} A \vee B$.

2.20. $\vdash_{ic} A \leftrightarrow B$ and $\vdash_{ic} B \leftrightarrow A$ if and only if A and B are syntactical equivalents in the classical implicative-conjunctive calculus.

2.21. $\vdash_{ic} (A \to (B \;\&\; C)) \leftrightarrow (A \to B) \;\&\; (A \to C)$.

2.22. $\vdash_{ic} ((A \;\&\; B) \to C) \leftrightarrow (A \to B \to C)$.

2.23. $\vdash_{ic} ((A_1 \;\&\; \ldots \;\&\; A_n) \;\&\; B) \leftrightarrow (A_1 \;\&\; \ldots \;\&\; A_n \;\&\; B)$.

2.24. $\vdash_{ic} B_1 \;\&\ldots\;\&\; B_n$ if and only if $\vdash_{ic} B_j$, for all $j$ such that $1 \le j \le n$.

2.25. $\vdash (C \vee (A \;\&\; B)) \leftrightarrow ((C \vee A) \;\&\; (C \vee B))$.

2.26. $\vdash ((A \;\&\; B) \vee C) \leftrightarrow ((A \vee C) \;\&\; (B \vee C))$.

2.27. Soundness Theorem for all the propositional calculi specified above (ST).

## 3. The classical implicative-disjunctive propositional calculus

We shall demonstrate some preparatory results with an eye on the intended completeness proof for the classical implicative-disjunctive calculus.

**3.1.** $(\Delta\,[V; B])^{B} \vdash_{id} (\Delta\,[V; A \to B])^{A \to B}$.

Proof.

(a) Let $\Delta$ [V; B] and $\Delta$ [V; $A \to B$] be both empty. In this case, 3.1 is $B \vdash_{id} A \to B$ and follows from Ax1 and DT.

(b) Let $\Delta$ [V; B] be empty and $\Delta$ [V; $A \to B$] be non-empty. In this case, 3.1 has the form $B \vdash_{id} C_1 \vee \ldots \vee C_n \vee (A \to B)$ and follows from $B \vdash_{id} (C_1 \vee \ldots \vee C_n) \vee (A \to B)$ (which is an instance of 2.14), 2.15 and RSE.

(c) Let $\Delta$ [V; B] and $\Delta$ [V; $A \to B$] be both non-empty. In this case, 3.1 has the form $B_1 \vee \ldots \vee B_n \vee B \vdash_{id} C_1 \vee \ldots \vee C_k \vee (A \to B)$, each of the $B_1,\ldots,B_n$ being one of the $C_1,\ldots,C_k$, by 2.1. By 2.16 and DT, $\vdash_{id} (B_1 \vee \ldots \vee B_n) \to (C_1 \vee \ldots \vee C_k)$. But $B \to A \to B$ is an instance of Ax1, and so, by 2.12, $(B_1 \vee \ldots \vee B_n) \vee B \vdash_{id} (C_1 \vee \ldots \vee C_k) \vee (A \to B)$. Hence 3.1, by 2.15 and RSE.

**3.2.** $(\Delta [V; A])^{\sim A} \vdash_{id} (\Delta [V; A \to B])^{(A \to B)}$, if V(A) = F.

Proof. Let V(A) = F. By 2.2, $\Delta$ [V; A] is non-empty and so, by 2.1, 3.2 has the form $A \to (A_1 \vee \ldots \vee A_n) \vdash_{id} C_1 \vee \ldots \vee C_k \vee (A \to B)$, each of the $A_1,\ldots,A_n$ being one of the $C_1,\ldots,C_k$. By 2.16 and DT, $\vdash_{id} (A_1 \vee \ldots \vee A_n) \to (C_1 \vee \ldots \vee C_k)$. Therefore, $A \to (A_1 \vee \ldots \vee A_n) \vdash_{id} A \to (C_1 \vee \ldots \vee C_k)$, by 2.6. But $A \to (C_1 \vee \ldots \vee C_k) \vdash_{id} (C_1 \vee \ldots \vee C_k) \vee (A \to B)$ is an instance of 2.13; hence $A \to (A_1 \vee \ldots \vee A_n) \vdash_{id} (C_1 \vee \ldots \vee C_k) \vee (A \to B)$, and so 3.2, by 2.15 and RSE.

**3.3.** $(\Delta [V; A])^{A}, (\Delta [V; B])^{\sim B} \vdash_{id} (\Delta [V; A \to B])^{\sim (A \to B)}$, if V(B) = F.

Proof. Let V(B) = F. In this case, $(\Delta [V; A \to B])^{\sim (A \to B)}$ has the form $(A \to B) \to (C_1 \vee \ldots \vee C_k)$, by 2.2. We have that $(\Delta [V; A])^{A}$ either is A or has the form $A_1 \vee \ldots \vee A_n \vee A$, each of the $A_1,\ldots,A_n$ being one of the $C_1,\ldots,C_k$, by 2.1; hence $(\Delta [V; A])^{A} \vdash_{id} C_1 \vee \ldots \vee C_k \vee A$, by 2.16. Therefore,

(i) $(\Delta [V; A])^{A} \vdash_{id} (C_1 \vee \ldots \vee C_k) \vee A$, by 2.15 and RSE.

By 2.2, $(\Delta [V; B])^{\sim B}$ has the form $B \to (B_1 \vee \ldots \vee B_j)$, each of the $B_1,\ldots,B_j$ being one of the $C_1,\ldots,C_k$,; consequently $\vdash_{id} (B_1 \vee \ldots \vee B_j) \to (C_1 \vee \ldots \vee C_k)$, by 2.16 and DT, and so

(ii) $(\Delta [V; B])^{\sim B} \vdash_{id} B \to (C_1 \vee \ldots \vee C_k)$, by 2.6. But

(iii) $(C_1 \vee \ldots \vee C_k) \vee A, B \rightarrow (C_1 \vee \ldots \vee C_k) \vdash_{id} (A \rightarrow B) \rightarrow (C_1 \vee \ldots \vee C_k)$ is an instance of 2.17.

Hence 3.3, by (i)-(iii).

**3.4.** $(\Delta\,[V; A])^{A} \vdash_{id} (\Delta\,[V; A \vee B)])^{A \vee B}$ and

$(\Delta\,[V; A])^{A} \vdash_{id} (\Delta\,[V; B \vee A)])^{B \vee A}$.

Proof. Straightforward, by 2.1 and 2.16.

**3.5.** $(\Delta\,[V; A])^{\sim A}, (\Delta\,[V; B])^{\sim B} \vdash_{id} (\Delta\,[V; A \vee B)])^{\sim(A \vee B)}$, if $V(A) = V(B) = F$.

Proof. Let $V(A) = V(B) = F$. In this case, by 2.2 and 2.1, $(\Delta\,[V; A \vee B)])^{\sim(A \vee B)}$ has the form $(A \vee B) \rightarrow (C_1 \vee \ldots \vee C_k)$, $(\Delta\,[V; A])^{\sim A}$ has the form $A \rightarrow (A_1 \vee \ldots \vee A_n)$, and $(\Delta\,[V; B])^{\sim B}$ has the form $B \rightarrow (B_1 \vee \ldots \vee B_j)$, each of the $A_1,\ldots, A_n, B_1,\ldots, B_j$ being one of the $C_1,\ldots, C_k$. By 2.16, $\vdash_{id} (A_1 \vee \ldots \vee A_n) \rightarrow (C_1 \vee \ldots \vee C_k)$ and $\vdash_{id} (B_1 \vee \ldots \vee B_j) \rightarrow (C_1 \vee \ldots \vee C_k)$. Therefore, $(\Delta\,[V; A])^{\sim A} \vdash_{id} A \rightarrow (C_1 \vee \ldots \vee C_k)$ and $(\Delta\,[V; B])^{\sim B} \vdash_{id} B \rightarrow (C_1 \vee \ldots \vee C_k)$, by 2.6 and DT; hence 3.5, by Ax6 and MP.

**3.6.** Let K be a finite set of implicative-disjunctive formulas; if $H \vdash_{id} (J)^{A}$, for all mutually exclusive sets H and J such that $H \cup J = K$, then $\vdash_{id} A$.

Proof. By induction on the finite cardinal *n* of K. Let K be a finite set of implicative-disjunctive formulas and

(i) $H \vdash_{id} (J)^{A}$, for all mutually exclusive sets H and J such that $H \cup J = K$.

(a) If $n = 0$, then, for all subsets H and J of K, $H = J = \varnothing$, and so $(J)^{A}$ is A; hence $\varnothing \vdash_{id} A$, that is to say, $\vdash_{id} A$, by (i).

(b) If $n > 0$, let $B_1,\ldots,B_n$ be the distinct elements of K arranged according to their relative position in the linear order R; let K* be the set of elements C of K such that $C \neq B_1$; and let H* and J* be any mutually exclusive sets such that $H^* \cup J^* = K^*$. It is easily verifiable that

(ii) $(H^* \cup \{B_1\}) \cup J^* = K$;

(iii) $H^* \cup \{B_1\}$ and J* are mutually exclusive;

(iv) $H^* \cup (J^* \cup \{B_1\}) = K$;

(v) $H^*$ and $J^* \cup \{B_1\}$ are mutually exclusive.

Thus, we may conclude that

(vi) $(H^* \cup \{B_1\}) \vdash_{id} (J^*)^A$, by (i)-(iii);

(vii) $H^* \vdash_{id} (J^* \cup \{B_1\})^A$, by (i), (iv) and (v).

Moreover,

(viii) $H^* \vdash_{id} B_1 \rightarrow (J^*)^A$, by (vi) and DT,

(ix) $(J^* \cup \{B_1\})^A \vdash_{id} B_1 \vee (J^*)^A$, by 2.16, and so

(x) $H^* \vdash_{id} B_1 \vee (J^*)^A$, by (vii) and (ix).

Therefore $H^* \vdash_{id} (J^*)^A$, by (viii), (x) and 2.18.

But we assumed $H^*$ and $J^*$ to be any mutually exclusive sets such that $H^* \cup J^* = K^*$ and the cardinal of $K^*$ to be $(n - 1)$; hence, by the induction hypothesis, $\vdash_{id} A$.

**3.7.** Let A be an implicative-disjunctive formula. If V(A) = T, then $\Gamma\,[V; A] \vdash_{id} (\Delta\,[V; A])^A$; if V(A) = F, then $\Gamma\,[V; A] \vdash_{id} (\Delta\,[V; A])^{\sim A}$.

Proof. By induction on the length of A.

(a) Let A be an atomic formula. If V(A) = T, then $\Gamma\,[V; A]$ is A and $(\Delta\,[V; A])^A$ is also A, and so $\Gamma\,[V; A] \vdash_{id} (\Delta\,[V; A])^A$. If V(A) = F, then $\Gamma\,[V; A]$ is empty and $(\Delta\,[V; A])^{\sim A}$ is $A \rightarrow A$; hence $\Gamma\,[V; A] \vdash_{id} (\Delta\,[V; A])^{\sim A}$, by 2.5.

(b) Let A be an implication $B \rightarrow C$ and V(C) = T. In this case, V(A) = T, $\Gamma\,[V, C] \vdash_{id} (\Delta\,[V, C])^C$, by the induction hypothesis, and $(\Delta\,[V, C])^C \vdash_{id} (\Delta\,[V, B \rightarrow C])^{B \rightarrow C}$, by 3.1. Since $\Gamma\,[V, C] \subseteq \Gamma\,[V, B \rightarrow C]$, it follows that $\Gamma\,[V, B \rightarrow C] \vdash_{id} (\Delta\,[V, B \rightarrow C])^{B \rightarrow C}$.

(c) Let A be an implication $B \rightarrow C$ and V(B) = F. In this case, V(A) = T, $\Gamma\,[V, B] \vdash_{id} (\Delta\,[V, B])^{\sim B}$, by the induction hypothesis, and $(\Delta\,[V, B])^{\sim B} \vdash_{id} (\Delta\,[V, B \rightarrow C])^{B \rightarrow C}$, by 3.2. Since $\Gamma\,[V, B] \subseteq \Gamma\,[V, B \rightarrow C]$, it follows that $\Gamma\,[V, B \rightarrow C] \vdash_{id} (\Delta\,[V, B \rightarrow C])^{B \rightarrow C}$.

(d) Let A be an implication $B \rightarrow C$, V(B) = T and V(C) = F. In this case, V(A) = F, $\Gamma\,[V, B] \vdash_{id} (\Delta\,[V, B])^B$ and $\Gamma\,[V, C] \vdash_{id} (\Delta\,[V, C])^{\sim C}$, by the induction hypothesis, and $(\Delta\ [V;\ B])^{B}$, $(\Delta\ [V;\ C])^{\sim C}$ $\vdash_{id}$ $(\Delta\ [V;\ B \rightarrow C])^{\sim(B \rightarrow C)}$, by 3.3. Since $\Gamma\ [V, B] \subseteq \Gamma\ [V, B \rightarrow C]$ and $\Gamma\ [V, C] \subseteq \Gamma\ [V, B \rightarrow C]$, by 2.1, it follows that $\Gamma\ [V, B \rightarrow C] \vdash_{id} (\Delta\ [V, B \rightarrow C])^{\sim(B \rightarrow C)}$.

(e) Let A be a disjunction $B \vee C$ and either V(B) = T or V(C) = T. In this case, V(A) = T, and we have that

(i) either $\Gamma\ [V, B] \vdash_{id} (\Delta\ [V, B])^{B}$, by the induction hypothesis, and $(\Delta\ [V, B])^{B} \vdash_{id} (\Delta\ [V, B \vee C])^{B \vee C}$, by 3.4, or $\Gamma\ [V, C] \vdash_{id} (\Delta\ [V, C])^{C}$, by the induction hypothesis, and $(\Delta\ [V, C])^{C} \vdash_{id} (\Delta\ [V, B \vee C])^{B \vee C}$, by 3.4.

Since $\Gamma\ [V, B] \subseteq \Gamma\ [V, B \rightarrow C]$ and $\Gamma\ [V, C] \subseteq \Gamma\ [V, B \rightarrow C]$, it follows from (i) that $\Gamma\ [V, B \vee C] \vdash_{id} (\Delta\ [V, B \vee C])^{B \vee C}$.

(f) Let A be a disjunction $B \vee C$ and V(B) = V(C) = F. In this case, V(A) = F, and we have that

(i) $\Gamma\ [V, B] \vdash_{id} (\Delta\ [V, B])^{\sim B}$ and $\Gamma\ [V, C] \vdash_{id} (\Delta\ [V, C])^{\sim C}$, by the induction hypothesis;
(ii) $(\Delta\ [V, B])^{\sim B}$, $(\Delta\ [V, C])^{\sim C} \vdash_{id} (\Delta\ [V, B \vee C])^{\sim(B \vee C)}$, by 3.5.

Hence $\Gamma\ [V, B]$, $\Gamma\ [V, C] \vdash_{id} (\Delta\ [V, B \vee C])^{\sim(B \vee C)}$, by (i) and (ii). Since $\Gamma\ [V, B] \subseteq \Gamma\ [V, B \vee C]$ and $\Gamma\ [V, C] \subseteq \Gamma\ [V, B \vee C]$, it follows that $\Gamma\ [V, B \vee C] \vdash_{id} (\Delta\ [V, B \vee C])^{\sim(B \vee C)}$.

**3.8.** Completeness theorem for the implicative-disjunctive calculus: for any implicative-disjunctive formula A and any finite set K of implicative disjunctive formulas, $\vdash_{id}$ A if and only if $\models$ A and K $\vdash_{id}$ A if and only if K $\models$ A.

Proof. Let A be an implicative-disjunctive formula.

(a) If $\vdash_{id}$ A, then $\models$ A, by ST.

(b) Let $\models$ A, let K be the set of atomic subformulas of A, and let H and J be any mutually exclusive sets such that $H \cup J = K$. Let V be a truth value assignment for L that fulfills the following condition: for all atomic subformulas B of A, V(B) = T if and only if B belongs to H, and so V(B) = F if and only if B belongs to J. Since V(A) = T, 3.7 ensures that $\Gamma\ [V; A] \vdash_{id} (\Delta\ [V; A])^{A.}$ But $\Gamma\ [V; A]$ is H and $(\Delta\ [V; A]$ is J; hence H $\vdash_{id} (J)^{A}$, from which $\vdash_{id}$ A follows, by 3.6.

From this and DT it follows that K $\vdash_{id}$ A if and only if K $\models$ A, for any implicative formula A and any finite set K of implicative formulas.

## 4. The classical positive propositional calculus

The completeness of the classical positive propositional calculus can be proved using the Kalmár-style proof procedure that we applied to the implicative-disjunctive calculus in Part 3 above. Indeed, it can be shown that propositions analogous to 3.1-3.9 hold for the classical positive propositional calculus. To do this, nothing is required but remembering that this calculus is an extension of the implicative-disjunctive calculus and showing that 4.1 and 4.2 below hold good.

**4.1.** $(\Delta [V; A])^{A}, (\Delta [V; B])^{B} \vdash (\Delta [V; A \,\&\, B])^{A \,\&\, B}$.

Proof.

(a) If $\Delta$ [V; A & B] is empty, then $(\Delta [V; A])^{A}$ and $(\Delta [V; B])^{B}$ are also empty. In this case, 4.1 is A, B $\vdash$ A & B and follows from Ax9, DT, and MP.

(b) If $\Delta$ [V; A & B] is non-empty, then $(\Delta [V; A \,\&\, B])^{A \,\&\, B}$ has the form $C_1 \vee \ldots \vee C_k \vee (A \,\&\, B)$; thus, $(\Delta [V; A])^{A} \vdash C_1 \vee \ldots \vee C_k \vee A$ and $(\Delta [V; B])^{B} \vdash C_1 \vee \ldots \vee C_k \vee B$, by 2.16. From this it follows that $(\Delta [V; A])^{A} \vdash (C_1 \vee \ldots \vee C_k) \vee A$ and $(\Delta [V; B])^{B} \vdash (C_1 \vee \ldots \vee C_k) \vee B$, by 2.15 and RSE. Hence, $(\Delta[V; A])^{A}, (\Delta [V; B])^{B} \vdash ((C_1 \vee \ldots \vee C_k) \vee A) \,\&\, ((C_1 \vee \ldots \vee C_k) \vee B)$, by Ax9 and MP, and so $(\Delta[V; A])^{A}, (\Delta [V; B])^{B} \vdash ((C_1 \vee \ldots \vee C_k) \vee (A \,\&\, B)$, by 2.25; therefore 4.1, by 2.15 and RSE.

**4.2.** If V(A) = F, then $(\Delta [V; A])^{\sim A} \vdash (\Delta [V; A \,\&\, B])^{\sim (A \,\&\, B)}$ and $(\Delta [V; A])^{\sim A} \vdash (\Delta [V; B \,\&\, A])^{\sim (B \,\&\, A)}$.

Proof. Let V(A) = F. In this case, by 2.1 and 2.2, $(\Delta [V; A])^{\sim A}$ has the form $A \rightarrow (A_1 \vee \ldots \vee A_n)$, $(\Delta [V; A \,\&\, B])^{\sim (A \,\&\, B)}$ has the form $(A \,\&\, B) \rightarrow (C_1 \vee \ldots \vee C_k)$ and $(\Delta [V; B \,\&\, A])^{\sim (B \,\&\, A)}$ is $(B \,\&\, A) \rightarrow (C_1 \vee \ldots \vee C_k)$, each of the $A_1,\ldots,A_n$ being one of the $C_1,\ldots,C_k$. But $\vdash (A_1 \vee \ldots \vee A_n) \rightarrow (C_1 \vee \ldots \vee C_k)$, by 2.16 and DT, and so

$(\Delta\ [V; A])^{\sim A} \vdash A \rightarrow (C_1 \vee \ldots \vee C_k)$, by 2.16, DT and 2.6. But $(A \& B) \rightarrow A$ and $(B \& A) \rightarrow A$ are instances of Ax7 and Ax8 respectively. Therefore 4.2, by 2.6.

The completeness of the classical positive propositional calculus can also be inferred from the completeness of the classical implicative-disjunctive propositional calculus by resorting to 4.3 and 4.4 below.

If A is a formula of L, $\gamma$ (A) will denote the formula that fulfills the following conditions:

(a) the connective & does not occur in $\gamma$ (A) inside any disjunction or conjunction;

(b) $\gamma$ (A) can be obtained from A by *n* successive applications of the following rules ($n \geq 0$):

(i) replace a subformula of the form $C \rightarrow (D \& E)$ by $(C \rightarrow D) \& (C \rightarrow E)$;

(ii) replace a subformula of the form $(C \& D) \rightarrow E$ by $C \rightarrow D \rightarrow E$;

(iii) replace a subformula of the form $C \vee (D \& E)$ by $(C \vee D) \& (C \vee E)$;

(iv) replace a subformula of the form $(C \& D) \vee E$ by $(C \vee E) \& (D \vee E)$.

For each formula A of L, the existence and unicity of $\gamma$(A) can be easily proved.

**4.3.** For each formula A of L, $\vdash A \leftrightarrow \gamma(A)$, and $\models A \leftrightarrow \gamma(A)$.

Proof. By induction on the length of A, through 2.21, 2.22, 2.25, 2.26, RSE and ST.

**4.4.** For each formula A of L, there are implicative-disjunctive formulas $B_1,\ldots,B_n$ such that $\vdash A \leftrightarrow (B_1 \& \ldots \& B_n)$ and $\models A \leftrightarrow (B_1 \& \ldots \& B_n)$.

Proof. Let A be any formula of L. Given 4.3 and RSE, it is enough to prove that

(i) there are implicative-disjunctive formulas $B_1,\ldots,B_n$ such that $\vdash \gamma(A) \leftrightarrow (B_1 \& \ldots \& B_n)$ and $\models \gamma(A) \leftrightarrow (B_1 \& \ldots \& B_n)$.

This can be done by induction on the length of $\gamma$(A).

(a) If $\gamma(A)$ is an atomic formula, an implication, or a disjunction, then $\gamma(A)$ is an implicative-disjunctive formula, by definition of $\gamma$, and (i) is trivially true, since any formula of L is syntactically and semantically equivalent to itself.

(b) If $\gamma(A)$ has the form C & D, then C is $\gamma(C)$ and D is $\gamma(D)$, by definition of $\gamma$. By the induction hypothesis, there are implicative-disjunctive formulas $C_1,\ldots,C_j$, $D_1,\ldots,D_k$ such that $\vdash C \leftrightarrow (C_1 \,\&\, \ldots \&\, C_j)$, $\models C \leftrightarrow (C_1 \,\&\, \ldots \&\, C_j)$, $\vdash D \leftrightarrow (D_1 \,\&\, \ldots \&\, D_k)$ and $\models D \leftrightarrow (D_1 \,\&\, \ldots \&\, D_k)$; Hence $\vdash \gamma(A) \leftrightarrow (C_1 \,\&\, \ldots \&\, C_j) \,\&\, (D_1 \,\&\, \ldots \&\, D_k)$ and $\models \gamma(A) \leftrightarrow (C_1 \,\&\, \ldots \&\, C_j) \,\&\, (D_1 \,\&\, \ldots \&\, D_k)$, by RSE. Consequently, $\vdash \gamma(A) \leftrightarrow (C_1 \,\&\, \ldots \&\, C_j \,\&\, D_1 \,\&\, \ldots \,\&\, D_k)$, by $(j-1)$ applications of 2.23 and RSE, and so $\models \gamma(A) \leftrightarrow (C_1 \,\&\, \ldots \&\, C_j \,\&\, D_1 \,\&\, \ldots \&\, D_k)$, by ST.

Now, ST ensures that all theses of the classical positive propositional calculus are tautologies. Let A be any tautology of L; by 4.4, there are implicative-disjunctive formulas $B_1,\ldots,B_n$ such that $\vdash A \leftrightarrow B_1 \,\&\ldots \,\&\, B_n$ and $\models A \leftrightarrow B_1 \,\&\ldots \,\&\, B_n$. Hence $\models B_1 \,\&\ldots \,\&\, B_n$ and it is obvious that, for all $j$ such that $1 \leq j \leq n$, $\models B_j$; by the completeness theorem for the implicative-disjunctive calculus, for all $j$ such that $1 \leq j \leq n$, $\vdash_{id} B_j$, and so $\vdash B_j$. Therefore $\vdash B_1 \,\&\ldots \,\&\, B_n$, by 2.24, and thus $\vdash A$, by 4.4 and RSE.

## 5. The classical implicative and the classical implicative-conjunctive propositional calculi.

The classical implicative-disjunctive propositional calculus can be mirrored in the classical implicative propositional calculus in a way that allows the completeness of the latter to be smoothly derived from the completeness of the former.

Let A be an implicative-disjunctive formula. We define the *implicative translation* of A, denoted by $\tau(A)$, as being the implicative formula that can be obtained from A by $n$ successive applications of the following rule ($n \geq 0$): replace a subformula of the form $B \vee C$ by $(B \rightarrow C) \rightarrow C$. For each formula A of L, the existence and unicity of $\tau(A)$ can be easily proved.

**5.1.** $(B \to D), (C \to D), (B \to C) \to C \vdash_i D$

Proof.

(i) $B \to D$, by hypothesis;

(ii) $C \to D$, by hypothesis;

(iii) $(B \to C) \to C$, by hypothesis;

(iv) $(D \to C) \to C$, by (i), (iii) and 2.10;

(v) $(D \to C) \to D$, by (ii), (iv), 2.6 and MP;

(vi) $((D \to C) \to D) \to D$ is an instance of Axb

(vii) D, by (v), (vi) and MP.

**5.2.** $\tau(A) \vdash_{id} A$ and $A \vdash_{id} \tau(A)$, for any implicative-disjunctive formula A.

Proof. By induction on the length of A, through 2.18, DT, 2.19, and RSE.

**5.3.** $\vdash_{id} A$ if and only if $\vdash_i \tau(A)$, for any implicative-disjunctive formula A.

Proof.

(a) If $\vdash_i \tau(A)$, then obviously $\vdash_{id} \tau(A)$, and so $\vdash_{id} A$, by 5.2.

(b) Now let A be any thesis of the classical implicative-disjunctive propositional calculus. We shall prove that $\vdash_i \tau(A)$ by induction on the length of a derivation of A in this calculus. If A is an instance of Ax1, Ax2, or Ax3, then $\tau(A)$ is also an instance of Ax1, Ax2, or Ax3. If A is an instance of Ax4 or Ax5, then $\tau(A)$ has either the form $\tau(B) \to (\tau(B) \to \tau(C)) \to \tau(C)$ or the form $\tau(B) \to (\tau(C) \to \tau(B)) \to \tau(B)$; hence $\vdash_i \tau(A)$, by 2.7 or 2.8. If A is an instance of Ax6, then $\tau(A)$ has the form $(\tau(B) \to \tau(D)) \to (\tau(C) \to \tau(D)) \to ((\tau(B) \to \tau(C)) \to \tau(C)) \to \tau(D)$, and thus $\vdash_i \tau(A)$, by 5.1 and DT. If A is derived from $B \to A$ and B, then $\vdash_i \tau(B \to A)$ and $\vdash_i \tau(B)$, by the induction hypothesis; but $\tau(B \to A)$ is $\tau(B) \to \tau(A)$; hence $\vdash_i \tau(A)$, by MP.

**5.4.** Completeness theorem for the implicative calculus: for any implicative formula A of L, $\vdash_i A$ if and only if $\models A$, and, for any finite set K of implicative formulas, $K \vdash_I A$ if and only if $K \models A$.

Proof. Let A be any implicative formula of L. If $\vdash_i$ A, then $\models$ A, by ST. Now, if $\models$ A, then $\vdash_{id}$ A, by the completeness theorem for the implicative-disjunctive calculus and so $\vdash_i \tau(A)$, by 5.3; but $\tau(A)$ is A itself, therefore $\vdash_i$ A. From this and DT it follows that K $\vdash_i$ B if and only if K $\models$ B, for any implicative formula B and any finite set K of implicative formulas.

Proposition 5.3 and RSE ensure that the completeness of the implicative propositional calculus can also be directly proved through a Kalmár-style procedure. This can be done by simply replacing in the proofs of 3.1-3.9 all references to the implicative-disjunctive propositional calculus by references to the implicative propositional calculus, and all references to implicative-disjunctive formulas by references to their respective implicative translations.

Finally, the completeness of the classical implicative-conjunctive propositional calculus can be inferred from the completeness of the classical implicative propositional calculus and 5.5 below exactly as the completeness of the classical positive propositional calculus was inferred from the completeness of the classical implicative-disjunctive propositional calculus and 4.4 in Part 4 above.

**5.5.** For any implicative-conjunctive formula A of L, there are implicative formulas $B_1,\ldots,B_n$ such that $\vdash_{ic} A \leftrightarrow (B_1 \,\&\ldots\&\, B_n)$ and $\models A \leftrightarrow (B_1 \,\&\ldots\&\, B_n)$.

Proof. Similar to the proof of 4.4.

**5.6.** Completeness theorem for the classical implicative-conjunctive propositional calculus: for any implicative-conjunctive formula of L, $\vdash_{ic}$ A if and only if $\models$ A, and, for any finite set K of implicative-conjunctive formulas, K $\vdash_{ic}$ A if and only if K $\models$ A.

Proof. We can infer 5.6 from 5.5 exactly as we inferred the completeness of the classical positive propositional calculus from 4.4 in Part 4 above.

It is worth noting that 4.4, 5.2, and RSE allow us to prove a stronger variant of 4.4.

**5.7.** For each formula A of L, there are implicative formulas $B_1,\ldots,B_n$ such that $\vdash A \leftrightarrow (B_1 \,\&\ldots\&\, B_n)$ and $\models A \leftrightarrow (B_1 \,\&\ldots\&\, B_n)$.